\documentclass[12pt]{article}
\usepackage{amsmath}
\usepackage{amssymb}
\usepackage{latexsym}
\usepackage{amsthm}
\usepackage{tabularx}
\usepackage{booktabs}
\usepackage{mathrsfs}
\usepackage{enumitem}
\usepackage{ytableau}
\usepackage{graphicx}
\usepackage{algorithm,algorithmic}
\usepackage[colorlinks,
linkcolor=red,
anchorcolor=blue,
citecolor=green]{hyperref}

\usepackage{epic}

\newtheorem{thm}{Theorem}[section]
\newtheorem{lem}[thm]{Lemma}
\newtheorem{prop}[thm]{Proposition}
\newtheorem{cor}[thm]{Corollary}
\newtheorem{conj}[thm]{Conjecture}

\newcommand{\p}{\overline{p}}

\newcommand{\M}{\overline{M}}

\newcommand{\bmu}{\overline{\mu}}

\allowdisplaybreaks

\begin{document}
	\begin{center}
		{\large \bf  The Briggs inequality for partitions and overpartitions}
	\end{center}
	
	\begin{center}
		Xin-Bei Liu$^1$, and Zhong-Xue Zhang$^{2}$\\[6pt]
		
		$^{1,2}$Center for Combinatorics, LPMC\\
		Nankai University, Tianjin 300071, P. R. China\\[6pt]

		Email: $^{1}${\tt liuxib@mail.nankai.edu.cn},
			   $^{2}${\tt zhzhx@mail.nankai.edu.cn}
		\end{center}
	\noindent\textbf{Abstract.}
A sequence of $\{a_n\}_{n\ge 0}$ satisfies the Briggs inequality if
\begin{align*}
a_n^2(a_n^2-a_{n-1}a_{n+1})>a_{n-1}^2(a_{n+1}^2-a_na_{n+2})
\end{align*}
holds for any $n\ge 1$. In this paper we show that
both the partition function $\{p(n+N_0)\}_{n\geq 0}$ and the overpartition function $\{\p(n+\overline{N}_0)\}_{n\ge 0}$ satisfy the Briggs inequality for some $N_0$ and $\overline{N}_{0}$.
Based on Chern's formula for $\eta$-quotients, we further prove that the $k$-regular partition function  $\{p_k(n+N_{k})\}_{n\geq 0}$ and the $k$-regular overpartition function $\{\p_k(n+\overline{N}_k)\}_{n\ge 0}$ also satisfy the Briggs inequality for $2\le k\le 9$ and some $N_k,\overline{N}_{k}$.

	\noindent \emph{AMS Classification 2020:} {05A20, 11B83}

	\noindent \emph{Keywords:} Briggs inequality, partition, overpartition, $k$-regular partition, $k$-regular overpartition
	
	\noindent \emph{The corresponding author:} Zhong-Xue Zhang 
	
	\noindent \emph{Suggested running title:} The Briggs inequality for partitions and overpartitions
	
\section{Introduction}

In the study of binding polynomials, Briggs \cite{Briggs-1985} proposed the following conjecture.
\begin{conj}
Suppose that $f(x)=\sum_{l=0}^{n}a_lx^l$ is a polynomial
with nonnegative coefficients.	
If $f(x)$ has only negative zeros, then, for $1\le l\le n-1$,
\begin{align}
	a_{l-1}a_{l+2}^2+a_l^2a_{l+3}+a_{l+1}^3&>a_{l+1}(a_{l-1}a_{l+3}+2a_{l}a_{l+2}),\label{2-log}\\[5pt]
	a_l^2(a_l^2-a_{l-1}a_{l+1})&>a_{l-1}^2(a_{l+1}^2-a_la_{l+2}).\label{Br-ineq}
\end{align}
\end{conj}

As pointed out by Zhang and Zhao \cite{Zhang-Zhao-2024},
the first inequality \eqref{2-log} can be derived from a result due to Br\"anden \cite{Branden-2015}. This paper is mainly concerned with the second inequality \eqref{Br-ineq}, which has been proved by Fan and Wang \cite{Fan-Wang} recently. Although Briggs' original conjecture is stated for a finite sequence, it is natural to explore such inequalities for infinite sequences. Following Zhang and Zhao \cite{Zhang-Zhao-2024}, we say that a sequence  $\{a_n\}_{n\ge 0}$ satisfies the Briggs inequality if
\begin{align*}
a_n^2(a_n^2-a_{n-1}a_{n+1})>a_{n-1}^2(a_{n+1}^2-a_na_{n+2})
\end{align*}
holds for any $n\ge 1$.
Zhang and Zhao \cite{Zhang-Zhao-2024} proved that both the Boros-Moll sequence and two of its variations satisfy the Briggs inequality.

In this paper, we aim to show that the Briggs inequality is satisfied by the partition function, the overpartition function, the $k$-regular partition function for $2\leq k\leq 9$, and the $k$-regular overpartition function for $2\leq k\leq 9$. Recall that a partition of $n$ is a weakly decreasing sequence of positive numbers whose sum is $n$. The partition function $p(n)$ counts the number of partitions of $n$. As a broad generalization of partitions, Corteel and Lovejoy \cite{Corteel-Lovejoy-2004} introduced the concept of overpartitions. By an overpartition of $n$ we mean a partition of $n$ such that the first occurrence of a number may be overlined. The overpartition function $\p(n)$ counts the number of overpartitions of $n$.
For example, there are three partitions of $3$, namely
$(3),(2,1),(1,1,1)$, and eight overpartitions of $3$, namely $(3),(\overline{3}),(2,1),(\overline{2},1),(2,\overline{1}),(\overline{2},\overline{1}),(1,1,1), (\overline{1},1,1)$.
Thus $p(3)=3$ and $\p(3)=8$.
For $k \geq 2$, a $k$-regular partition of $n$ is a partition of $n$ with no part divisible by $k$. A $k$-regular overpartition of $n$ can be defined in the same manner. As usual, we use $p_k(n)$ and $\p_k(n)$ to represent the $k$-regular partition function and the $k$-regular overpartition function, respectively.

Various interesting inequalities have been established for the partition function, the overpartition function, the $k$-regular partition function and the $k$-regular overpartition function. The log-concavity of  $\{p(n)\}_{n\ge 26}$  was independently proved by Nicolas in \cite{Nicolas-1978} and by DeSalvo and Pak in \cite{DeSalvo-Pak-2015}. Chen, Jia and Wang \cite{Chen-Jia-Wang-2019} showed that $\{p(n)\}_{n\ge 95}$ also satisfies higher order Tur\'an inequalities. For more information on
higher order Tur\'an inequalities and double Tur\'an inequalities, see \cite{Csordas-2003,Dimitrov-1998,Niculescu-2000,Rosset-1989}.
Furthermore, Chen, Jia and Wang conjectured that for $d\ge 4$ there is a positive number $N_d$ such that the order $d$ Tur\'an inequalities are valid for $p(n)$ when $n\geq N_d$.
Later, this conjecture was solved by Griffin, Ono, Rolen, and Zagier \cite{G-O-R-Z-2019}.
Hou and Zhang \cite{Hou-Zhang-2019} proved the asymptotic $r$-log-concavity of $\{p(n)\}_{n\geq 1}$ for any $r\geq 1$.
In particular, they established the $2$-log-concavity of $\{p(n)\}_{n\ge 222}$, which was independently proved by
Jia and Wang \cite{Jia-Wang-2020}.
For the overpartition function, Engel \cite{Engel-2017} proved that $\{\p(n)\}_{n\ge 2}$ is log-concave.
Liu and Zhang \cite{Liu-Zhang-2021} showed that the higher order Tur\'an inequalities are satisfied by $\{\p(n)\}_{n\ge 16}$. Following the work in \cite{Jia-Wang-2020}, Mukherjee \cite{Mukherjee-2023} showed that
$\{\p(n)\}_{n\ge 42}$ satisfies the double Tur\'an inequalities. Later,
Mukherjee, Zhang and Zhong \cite{Mukherjee-Zhang-Zhong-2023} proved the asymptotic $r$-log-concavity of
$\{\p(n)\}_{n\ge 1}$. By employing the result in \cite{G-O-R-Z-2019}, Craig
and Pun \cite{Craig-Pun-2021} showed that for $k\geq 2$, the $k$-regular partition function $p_k(n)$ satisfies the order $d$ Tur\'an inequalities for sufficiently large $n$.
{Furthermore, they conjectured that $\{p_2(n)\}_{n\ge 33}$ is log-concave and $\{p_2(n)\}_{n\ge 121}$ satisfies
the higher order Tur\'an inequalities.}
Based on Chern's asymptotic formula \cite{Chern-2019}, Dong and Ji \cite{Dong-Ji-2024} showed that  $\{p_k(n)\}_{n\ge N_k}$ is log-concave and satisfies higher order inequalities for $2\leq k\leq 5$ and some $N_k$, thus confirming the conjectures of Craig and Pun. Wang and Yang \cite{Wang-Yang} showed that  $\{p_2(n)\}_{n\ge 271}$ satisfies double Tur\'an inequalities, as conjectured by Dong and Ji \cite{Dong-Ji-2024}.
Peng, Zhang and Zhong \cite{Peng-Zhang-Zhong-2023} proved that  $\{\p_k(n)\}_{n\ge \overline{N}_k}$ is log-concave and satisfies higher order inequalities for $2\leq k\leq 9$ and some $\overline{N}_k$.
	
We would like to point out that, for a given positive sequence $\{a_n\}_{n\ge 0}$, the inequality \eqref{2-log} is equivalent to the double Tur\'an inequality,
 while the Briggs inequality \eqref{Br-ineq}
is closely related to the log-concavity.
Recall that a sequence $\{a_n\}_{n\ge 0}$ is said to be log-concave if $a_{n}^2-a_{n+1}a_{n-2}\ge 0$ for any $n\ge 1$. Note that if a log-concave sequence $\{a_n\}_{n\ge 0}$ satisfies
	\begin{align}\label{main}
		a_{n+1}(a_n^2-a_{n-1}a_{n+1})>a_{n-1}(a_{n+1}^2-a_na_{n+2}),
	\end{align}
then it also satisfies the Briggs inequality. 
This follows from the fact that the log-concavity of $\{a_n\}_{n\ge 0}$ implies $a_n^2 \geq a_{n+1}a_{n-1}$, and therefore
$$a_n^2(a_n^2-a_{n-1}a_{n+1})\geq a_{n+1}a_{n-1}(a_n^2-a_{n-1}a_{n+1})>a_{n-1}^2(a_{n+1}^2-a_na_{n+2}).$$
In order to show that
$p(n), \p(n), p_k(n)$ and $\p_k(n)$ (after ignoring some initial terms) satisfy the Briggs inequality, it suffices to show that they satisfy the stronger inequality \eqref{main} in view of the aforementioned log-concavity of these partition functions.

The remainder of this paper is organized as follows.
In Section \ref{sec-partition} we prove the Briggs inequality of the partition function and the overpartition function by using the bounds of $p(n)$ and $\p(n)$ given by Wang and Yang \cite{Wang-Yang-2023}.
In Section \ref{sec-regular partition} we show that, for $2\leq k\leq 9$,
the $k$-regular partition function $p_k(n)$ and the $k$-regular overpartition function $\p_k(n)$
satisfy the Briggs inequality
by using some explicit bounds of $p_k(n)$ and $\p_k(n)$,
which can be obtained from Chern's formula of $\eta$-quotients \cite{Chern-2019}.

\section{Partition functions}\label{sec-partition}

The main objective of this section is to prove
that both the partition function and the overpartition function satisfy the Briggs inequality. Taking $a_n$ to be the partition function $p(n)$ or $\p(n)$, we only need to prove inequality \eqref{main}, as discussed earlier.
Note that \eqref{main} can rewritten as
		\begin{align}\label{main-equivalent-form}
a_{n+1}a_n^2-2a_{n-1}a_{n+1}^2+a_{n-1}a_na_{n+2}>0.
		\end{align}

\subsection{Partitions}

For the partition function $p(n)$, we have the following result.

	\begin{thm}\label{thm2}
	For all $n\geq 114$, we have
		\begin{align}\label{eq-main-par}
			p(n+1)p(n)^2-2p(n-1)p(n+1)^2+p(n-1)p(n)p(n+2)>0.
		\end{align}
	\end{thm}

To prove Theorem \ref{thm2}, we use the following
upper and lower bounds for $p(n)$ given by Wang and Yang \cite{Wang-Yang-2023}.
Let $\mu(n)=\dfrac{\pi\sqrt{24n-1}}{6}$ and
	\begin{align}\label{eq-f-g}
		f(t)=\frac{1}{t^2}\left( 1-\frac{1}{t}-\frac{1}{t^{10}}\right),\quad g(t)=\frac{1}{t^2}\left( 1-\frac{1}{t}+\frac{1}{t^{10}}\right).
	\end{align}
Wang and Yang \cite{Wang-Yang-2023} obtained the following result.

 	\begin{lem}[\cite{Wang-Yang-2023}, Lemma 2.1]\label{lem-WY-par}
 Let $\mu(n),f(t),g(t)$ be defined as above. Then for all  $n\geq 1520$, we have
		\begin{align*}
			\dfrac{\sqrt{12}\pi^2e^{\mu(n)}}{36}f\left(\mu(n)\right)<p(n)<	\dfrac{\sqrt{12}\pi^2e^{\mu(n)}}{36}g\left(\mu(n)\right).
		\end{align*}
	\end{lem}

Note that $p(n-1), p(n), p(n+1)$ and $p(n+2)$ appear in \eqref{eq-main-par}. 
In order to use the bounds of these values given by Lemma \ref{lem-WY-par}, for notational convenience, we set
	\begin{align}\label{eq-xi}
		x=\mu(n-1),\quad
		x_1=\mu(n),\quad
		x_2=\mu(n+1),\quad x_3=\mu(n+2)
	\end{align}
throughout this subsection.
In our proof of \eqref{eq-main-par}, the value of $x$ will be used to estimate $x_1,x_2$ and $x_3$. We have the following result.

\begin{lem}\label{lem-x-estimation}
Let $x, x_1,x_2$ and $x_3$ be functions of $n$ as defined in \eqref{eq-xi}, and let
		\begin{align}
			x_{11}=\check{h}_x\left(\frac{2\pi^2}{3}\right),\quad x_{21}&=\check{h}_x\left(\frac{4\pi^2}{3}\right),\quad x_{31}=\check{h}_x\left(2\pi^2\right),\label{eq-xi1-formula}\\
			x_{12}=\hat{h}_x\left(\frac{2\pi^2}{3}\right),\quad x_{22}&=\hat{h}_x\left(\frac{4\pi^2}{3}\right),\quad x_{32}=\hat{h}_x\left(2\pi^2\right),\label{eq-xi2-formula}
		\end{align}
where
 	\begin{align}
		&\check{h}_x(a):=
		x+\frac{a}{2x}-\frac{a^2}{8x^3}+\frac{a^3}{16x^5}-\frac{5a^4}{64x^7}\label{eq-x-},\\
		&\hat{h}_x(a):=x+\frac{a}{2x}-\frac{a^2}{8x^3}+\frac{a^3}{16x^5}\label{eq-x+}.	
	\end{align}
Then, for $x>\sqrt{2}\pi$ and hence $n\geq 5$,
		\begin{align*}
			x_{i1}<x_{i}<x_{i2}
		\end{align*}
holds for $1\leq i\leq 3$.
\end{lem}

\begin{proof}
One can directly verify that
	\begin{align*}
		x_1=\sqrt{x^2+\frac{2\pi^2}{3}},\quad  x_2=\sqrt{x^2+\frac{4\pi^2}{3}},\quad x_3=\sqrt{x^2+2\pi^2}.
	\end{align*}
Since each of $x_i$ is of the form $\sqrt{x^2+a}$ for some positive number $a$, it suffices to show that
$\check{h}_x(a)<\sqrt{x^2+a}<\hat{h}_x(a)$
for $x^2>a$.
Keeping in mind that $x$ is always positive, Newton's binomial theorem tells that
	\begin{align*}
\sqrt{x^2+a}&=x\left(1+\frac{a}{x^2}\right)^{\frac{1}{2}}=x\left(\sum_{k\geq 0} \binom{\frac{1}{2}}{k}\left(\frac{a}{x^2}\right)^k\right)\\[5pt]
&=x+\frac{a}{2x}-\frac{a^2}{8x^3}+\frac{a^3}{16x^5}-\frac{5a^4}{128x^7}+\frac{7a^5}{256x^9}-\frac{21a^6}{1024x^{11}}+O\left(\frac{1}{x^{13}}\right).
	\end{align*}
By considering the difference of two adjacent terms in the above expansion, one can show that if $x^2>a$ then $\check{h}_x(a)<\sqrt{x^2+a}<\hat{h}_x(a)$. This completes the proof.
\end{proof}

Now we are in the position to prove Theorem \ref{thm2}.

\begin{proof}[Proof of Theorem \ref{thm2}]
For $114\leq n\leq 1520$ one can directly verify \eqref{eq-main-par}. From now on we assume that $n\geq 1521$, or equivalently, $x\geq 100$.

Recalling the bounds of $p(n)$ given in Lemma \ref{lem-WY-par}, we obtain that
\begin{align*}
p(n+1)p(n)^2-2p(n-1)p(n+1)^2+p(n-1)p(n)p(n+2)>\left( \frac{\sqrt{12}\pi^2}{36}\right)^3 F_1(x),
		\end{align*}
		where
		\begin{align}\label{eq-F-1}
F_1(x)=e^{x_{21}+2x_{11}}f(x_2)f(x_{1})^2-2e^{2x_{22}+x}g(x_2)^2g(x)+e^{x+x_{11}+x_{31}}f(x)f(x_1)f(x_3),
		\end{align}
the symbols $x,x_1,x_2,x_3,x_{11},x_{21},x_{31},x_{22}$ are defined as in \eqref{eq-xi}, \eqref{eq-xi1-formula} and  \eqref{eq-xi2-formula}, and the functions $f(t),g(t)$ are given by \eqref{eq-f-g}.

It remains to show that $F_1(x)>0$.
Let $z_2=\hat{h}_x(\frac{8\pi^2}{9})$.
Note that the sign of $F_1(x)$ coincides with
that of $F_1(x)e^{-3z_2}$. We find that it is more convenient to deal with the latter.
It is routine to verify that, for $x\geq 4$,
		\begin{align*}
x_{21}+2x_{11}-3z_2&=\frac{-\pi^4(18x^4-26\pi^2x^2+135\pi^4)}{486x^7}<0,\\[5pt]
			2x_{22}+x-3z_2&=\frac{-4\pi^4 \left(9x^2-10\pi^2\right)}{243 x^5}<0,\\[5pt]
x+x_{11}+x_{31}-3z_2&=\frac{-\pi^4(126x^4-188\pi^2x^2+615\pi^4)}{486 x^7}<0.
		\end{align*}
While, for $t<0$, we have
	\begin{align} \label{eq-Ebound}
		G_1(t)<e^t<G_2(t),
	\end{align}
	where
	\begin{align}\label{eq-e1-e2}	
		G_1(t)=1+t+\frac{t^2}{2}+\frac{t^3}{6},\quad G_2(t)=1+t+\frac{t^2}{2}.
	\end{align}
Thus, for $x\geq 4$, there holds
		\begin{align}
			&e^{x_{21}+2x_{11}-3z_2}> G_1(x_{21}+2x_{11}-3z_2),\notag\\
			&e^{2x_{22}+x-3z_2}<G_2(2x_{22}+x-3z_2),\label{eq-ex}\\
			&e^{x+x_{11}+x_{31}-3z_2}> G_1(x+x_{11}+x_{31}-3z_2).\notag
		\end{align}

To prove that $F_1(x)e^{-3z_2}>0$, we also need to estimate the values of $f(x_i)$ and $g(x_i)$.
For $i=1,2,3$ letting
		\begin{align*}
			f_{i}(t)=\dfrac{t^{10}-t^8x_{i 2}-1}{t^{12}}\qquad\mbox{and}\qquad g_{i}(t)=\dfrac{t^{10}-t^8x_{i 1}+1}{t^{12}},
		\end{align*}
one can check that for $x\ge 1$,
		\begin{align}\label{eq-fbound}
			f_{i}(x_i)<f(x_i)\qquad \mbox{and}\qquad g(x_i)<g_{i}(x_i).
		\end{align}

Combining \eqref{eq-F-1}, \eqref{eq-ex}
and \eqref{eq-fbound}, we see that
		\begin{align*}
			F_1(x)e^{-3z_2}
			>F_2(x)
		\end{align*}
holds for $x\geq 4$,  where
		\begin{align*}
F_2(x)=&\big(G_1(x_{21}+2x_{11}-3z_2)f_{2}(x_2)f_1(x_{1})^2-2G_2(2x_{22}+x-3z_2)g_{2}(x_2)^2g(x)\\
			&\quad+G_1(x+x_{11}+x_{31}-3z_2)f(x)f_1(x_1)f_{3}(x_3)\big).
		\end{align*}
It turns out that $F_2(x)$, after simplification, can be written as the following form:
		\begin{align*}
			F_2(x)=\dfrac{\sum_{l=0}^{92}a_lx^l}{2^{5}3^{19}x^{43}(x^2+2\pi^2)^6(3x^2+2\pi^2)^{12}(3x^2+4\pi^2)^{12}},
		\end{align*}
		where $a_l$ are the known numbers, and the values of $a_{92}$, $a_{91}$, $a_{90}$ are given below
		\begin{align*}
			a_{92}=2^{5}3^{41}\pi^6,\quad a_{91}=-2^{5}3^{37}(675\pi^6+7\pi^8),\quad a_{90}=2^{7}3^{38}(162\pi^6+127\pi^8).
		\end{align*}
With the help of mathematical software, we find that the largest real zero of $F_2(x)$ is less than $15$. Thus, for $x\geq 100$ we have $F_2(x)>0$. This completes the proof.
	\end{proof}

Based on Theorem \ref{thm2} we obtain the following result.

\begin{cor}\label{cor-par}
The sequence $\{p(n)\}_{n\ge 114}$ satisfies the Briggs inequality.
\end{cor}

\begin{proof}
From Theorem \ref{thm2} it follows that
\begin{align*}
p(n+1)(p(n)^2-p(n-1)p(n+1))&>p(n-1)(p(n+1)^2-p(n)p(n+2))
\end{align*}
for $n\geq 114$. By the log-concavity of  $\{p(n)\}_{n\ge 26}$ proved in \cite{Nicolas-1978} and  \cite{DeSalvo-Pak-2015}, we immediately obtain the desired result.
\end{proof}

\subsection{Overpartitions}

For the overpartition function $\p(n)$ we have the following result.

	\begin{thm}\label{thm-overpar}
	For all $n\geq 18$, we have
	\begin{align}\label{eq-main-overpar}
		\p(n+1)\p(n)^2-2\p(n-1)\p(n+1)^2+\p(n-1)\p(n)\p(n+2)>0.
	\end{align}
\end{thm}

It's worth noting that our approach to the partition function can be carried over verbatim to the overpartition function. Along the lines of the proof of Theorem \ref{thm2}, we first recall the upper and lower bounds for $\p(n)$ given by Wang and Yang \cite{Wang-Yang-2023}.

	\begin{lem}[\cite{Wang-Yang-2023}, Lemma 5.1]\label{lem-WY-over}
Let $\overline{\mu}(n)=\pi\sqrt{n}$. Then for any $n\geq 821$, we have
		\begin{align*}
			\dfrac{\pi^2e^{\overline{\mu}(n)}}{8}f(\bmu(n)) <\p(n)<\dfrac{\pi^2e^{\bmu(n)}}{8}g(\bmu(n)).
		\end{align*}
	\end{lem}

Now the proof of Theorem \ref{thm-overpar} can be given in the same manner as that of Theorem \ref{thm2}.

\begin{proof}[Proof of Theorem \ref{thm-overpar}]
Set
		$
		\overline{x}=\bmu(n-1),
		\overline{x}_1=\bmu(n),
		\overline{x}_2=\bmu(n+1),
		\overline{x}_3=\bmu(n+2).
		$
		Then
		\begin{align*}	
			\overline{x}_1=\sqrt{\overline{x}^2+\pi^2},\quad  \overline{x}_2=\sqrt{\overline{x}^2+2\pi^2},\quad \overline{x}_3=\sqrt{\overline{x}^2+3\pi^2}.
		\end{align*}
One can directly check that \eqref{eq-main-overpar} holds for $18\leq n\leq 821$.
Now we may assume that $n\geq 822$, or equivalently, $\overline x\geq 90$.

Using the functions given by \eqref{eq-x-} and \eqref{eq-x+}, we set
\begin{align*}
			\overline{x}_{11}=\check{h}_{\overline{x}}(\pi^2),\quad \overline{x}_{21} &=\check{h}_{\overline{x}}(2\pi^2),\quad \overline{x}_{31}=\check{h}_{\overline{x}}(3\pi^2),\\
			\overline{x}_{12}=\hat{h}_{\overline{x}}(\pi^2),\quad \overline{x}_{22} &=\hat{h}_{\overline{x}}(2\pi^2),\quad \overline{x}_{32}=\hat{h}_{\overline{x}}(3\pi^2).
		\end{align*}
		and
		\begin{align*}
\overline{f}_{i}(t)=\dfrac{t^{10}-t^8\overline{x}_{i2}-1}{t^{12}}, \quad \overline{g}_{i}(t)=\dfrac{t^{10}-t^8\overline{x}_{i1}+1}{t^{12}}, \quad \overline{z}_2=\hat{h}_{\overline{x}}(\frac{4\pi^2}{3}).
		\end{align*}

Similar to the proof of Lemma \ref{lem-x-estimation}, we can show that if $\overline x\geq 6$ then $\overline{x}_{i1}<\overline{x}_{i}<\overline{x}_{i2}$
for $1\leq i\leq 3$. One can also show that if $\overline{x}\ge 1$ then  $\overline{f}_{i}(\overline{x}_i)<f(\overline{x}_i)$ and $g(\overline{x}_i)<\overline{g}_{i}(\overline{x}_i)$,
where $f(t)$ and $g(t)$ are defined by \eqref{eq-f-g}.
A little computation shows that $\overline{x}_{21}+2\overline{x}_{11}-3\overline{z}_2<0,\, 2\overline{x}_{22}+\overline{x}-3\overline{z}_2<0$ and $x+\overline{x}_{11}+\overline{x}_{31}-3\overline{z}_2<0$ for $\overline x\geq 5$.

Let
			\begin{align*}
				\overline{F}_2(\overline{x})=&\big(G_1(\overline{x}_{21}+2\overline{x}_{11}-3\overline{z}_2)\overline{f}_{2}(\overline{x}_2)\overline{f}_1(\overline{x}_{1})^2-2G_2(2\overline{x}_{22}+\overline{x}-3\overline{z}_2)\overline{g}_{2}(\overline{x}_2)^2\overline{g}(\overline{x})\\
				&\quad+G_1(\overline{x}+\overline{x}_{11}+\overline{x}_{31}-3\overline{z}_2){f}(\overline{x})\overline{f}_1(\overline{x}_1)\overline{f}_{3}(\overline{x}_3)\big),
			\end{align*}
where $G_1(t)$ and $G_2(t)$ are defined as in \eqref{eq-e1-e2}.

By using the same arguments as in the proof of Theorem \ref{thm2}, we find that the inequality \eqref{eq-main-overpar} is equivalent to the positivity of
$\overline{F}_2(\overline{x})$.
After some simplification, we see that
			\begin{align*}
				\overline{F}_2(\overline{x})=\dfrac{\sum_{l=0}^{92}a_l\overline{x}^l}{2^{25}3^{7}x^{43} \left(\overline{x}^2+\pi ^2\right)^{12} \left(\overline{x}^2+2 \pi ^2\right)^{12} \left(\overline{x}^2+3 \pi ^2\right)^6},
			\end{align*}
			where $a_l$ are the known number, and the values of $a_{92}$, $a_{91}$, $a_{90}$ are given below
			\begin{align*}
				a_{92}=2^{22}3^{8}\pi^6,\quad a_{91}=-2^{21}3^{5}(450\pi^6+7\pi^8),\quad a_{90}=2^{23}3^{6}(108\pi^6+127\pi^8).
			\end{align*}
By using mathematical software, one can check that the largest real zero of $\overline{F}_2(\overline{x})$ is less than $18$.
Thus $\overline{F}_2(\overline{x})>0$ for $\overline x\geq 90$, as desired. This completes the proof.		
	\end{proof}

Based on Theorem \ref{thm-overpar} and the log-concavity of $\{\p(n)\}_{n\ge 2}$ due to Engel \cite{Engel-2017}, we immediately obtain the following result.

\begin{cor}
The sequence $\{\p(n)\}_{n\ge 18}$ satisfies the Briggs inequality.
\end{cor}
		
\section{$k$-regular partition functions }\label{sec-regular partition}

In this section we aim to prove that, for $2\leq k\leq 9$,
both the $k$-regular partition and the $k$-regular overpartition satisfy the Briggs inequality. Here our approach is exactly the same as that for the partition function and the overpartition given in Section \ref{sec-partition}.

Fixing an integer $2\leq k\leq 9$, let $a_n$ be either $p_k(n)$ or $\overline{p}_k(n)$.
In order to prove that the sequence  $\{a_n\}_{n\ge 0}$
satisfies the Briggs inequality, it suffices to show that it is log-concave and moreover it satisfies \eqref{main-equivalent-form}.
Dong and Ji \cite{Dong-Ji-2024} proved the log-concavity of ${p}_k(n)$ for $k=2$, and their proof also works for $3\leq k\leq 9$, provided that the appropriate upper and lower bounds of ${p}_k(n)$ are established.
The log-concavity of $\overline{p}_k(n)$ has been established by Peng, Zhang and Zhong  \cite{Peng-Zhang-Zhong-2023} when $2\leq k\leq 9$.
The proofs of Theorems \ref{thm2} and \ref{thm-overpar}
reveal that, to prove \eqref{main-equivalent-form} for
$a_n=p_k(n)$ and $a_n=\overline{p}_k(n)$, it is necessary to establish some bounds of $p_k(n)$ and $\overline{p}_k(n)$.
Dong and Ji \cite{Dong-Ji-2024} gave
certain upper and lower bounds of ${p}_2(n)$.
For our purpose, we shall present some explicit bounds of ${p}_k(n)$ when $3\leq k\leq 9$.
For the $k$-regular overpartitions where $2\leq k\leq 9$, we shall use the bounds of $\overline{p}_k(n)$ given by Peng, Zhang and Zhong \cite{Peng-Zhang-Zhong-2023}.


\subsection{$k$-regular partitions}\label{sec-k-par}
For the $k$-regular partition function $p(n)$, we obtain the following result.

	\begin{thm}\label{thm-p_kmain}
		For $2\leq k\leq 9$ and $n\geq N_k$, we have
		\begin{align}\label{p_kmain}
			p_k(n+1)p_k(n)^2-2p_k(n-1)p_k(n+1)^2+p_k(n-1)p_k(n)p_k(n+2)> 0,
		\end{align}
		where
		\begin{align*}
			&N_2=150,\quad  N_3=220,\quad N_4=75,\quad N_5=164,\quad\\
			& N_6=60,\quad N_7=148,\quad N_8=78,\quad N_9=138.
		\end{align*}
	\end{thm}

To prove the above theorm, we need to use the following explicit bounds of ${p}_k(n)$.
To maintain readability, we will provide the proofs for the explicit bounds in the appendix, as they are quite tedious.

	\begin{thm}\label{thm-p_k2}
 For $2\le k\le 9$, let
 \begin{align} \label{eq-ykmkn}
  \mu_k(n)=\frac{\pi}{6}\sqrt{\left(1-\frac{1}{k}\right)(24n+k-1)}, \quad M_k(n)=\frac{(k-1)\pi^2}{3k\sqrt{k}\mu_k(n)} I_1(\mu_k(n)),
  \end{align}
 where $I_1(s)$ denotes the first modified Bessel function of the first kind.
Then
	\begin{align}\label{eq-p-k-upper-lower}
			M_k(n)\left(1-\frac{1}{\mu_k(n)^6}\right)<p_k(n)<M_k(n)\left(1+\frac{1}{\mu_k(n)^6}\right),
		\end{align}
whenever $n\ge n_k$, where
	\begin{align*}
		&n_2=1067,\quad n_3=821,\quad n_4=711,\quad n_5=695,\quad\\
		&n_6=677,\quad n_7=652,\quad n_8=651,\quad n_9=615.
	\end{align*}
	\end{thm}

We also need the following bounds of the first modified Bessel function of the first kind, due to Dong and Ji \cite{Dong-Ji-2024}.

\begin{lem}[\cite{Dong-Ji-2024}, Lemma 2.2 and equation(3.14)]\label{Ji-Dong}
	Let
$I_1(s)$ denote the first modified Bessel function of the first kind, and let
	\begin{equation}\label{defi-I-e}
		D_I(s):=1-\frac{3}{8s}-\frac{15}{128s^2}-\frac{105}{1024s^3}
		-\frac{4725}{32768s^4}-\frac{72765}{262144s^5}.
	\end{equation}
	Then for $s\geq 26$, we have
	\begin{align*}
		\frac{e^s}{\sqrt{2\pi s}}\left( D_I(s)-\frac{31}{s^6}\right)	\leq	I_1(s) \leq \frac{e^s}{\sqrt{2\pi s}}\left( D_I(s)+\frac{31}{s^6}\right).
	\end{align*}
	Moreover,
	\begin{align}\label{eq-I1s}
		I_1(s)\geq \frac{e^s}{\sqrt{2\pi s}}\left(1-\frac{1}{2s}\right).
	\end{align}
\end{lem}

For notational convenience, let
\begin{align}\label{eq-y_k}
	x=\mu_6(n-1),\quad y=\mu_6(n),\quad
	z=\mu_6(n+1),\quad
	w=\mu_6(n+2)
\end{align}
throughout this subsection, where
$\mu_6(n)$ is given by \eqref{eq-ykmkn}.
Before giving the proof of Theorem \ref{thm-p_kmain}, let us estimate the values $x,z$ and $w$ in terms of $y$.
By \eqref{eq-ykmkn} we have
	\begin{align*}
		x=\sqrt{y^2-\frac{5\pi^2}{9}},\quad  z=\sqrt{y^2+\frac{5\pi^2}{9}},\quad w=\sqrt{y^2+\frac{10\pi^2}{9}}.
	\end{align*}
The following result is analogous to Lemma \ref{lem-x-estimation}, and its proof is omitted here.

\begin{lem}
    Let $x,y,z$ and $w$ be defined as in \eqref{eq-y_k},
    let $\check{h}_{y}(a)$ and $\hat{h}_{y}(a)$ be given by \eqref{eq-x-} and \eqref{eq-x+} respectively, and let
	\begin{align}\label{+--bound}
		x_1=\check{h}_{y}\left(-\frac{5\pi}{9}\right),\quad
		z_1=\check{h}_{y}\left(\frac{5\pi}{9}\right),\quad
		w_1=\check{h}_{y}\left(\frac{10\pi}{9}\right),\notag\\
		x_2=\hat{h}_{y}\left(-\frac{5\pi}{9}\right),\quad
		z_2=\hat{h}_{y}\left(\frac{5\pi}{9}\right),\quad
		w_2=\hat{h}_{y}\left(\frac{10\pi}{9}\right).
	\end{align}
Then, for $y\geq 3$, we have
	\begin{align}\label{y bound}
		x_1<x<x_2,\quad
		z_1<z<z_2,\quad
		w_1<w<w_2.
	\end{align}
\end{lem}

We proceed to give a proof of Theorem \ref{thm-p_kmain}.

\begin{proof}[Proof of Theorem \ref{thm-p_kmain}]
We shall take $k=6$ as an example to illustrate our proof, and the proofs for other values of $k$ can be carried out in the same manner.
Let $x, z$ and $w$ be defined in \eqref{eq-y_k}.
By Theorem \ref{thm-p_k2}, we find that for $n\ge 679$, i.e., $y\ge 61$,
			$$\frac{5\pi^2}{18\sqrt{6}y}I_1(y)\left(1-\frac{1}{y^6}\right)<p_6(n)<\frac{5\pi^2}{18\sqrt{6}y}I_1(y)\left(1+\frac{1}{y^6}\right).$$
By further applying Lemma \ref{Ji-Dong}, we get that
		\begin{align}\label{p6-bound}
			\frac{5\pi^{\frac{3}{2}}\cdot e^{y}}{36\sqrt{3}y^{\frac{3}{2}}}f(y)
			<p_6(n)<
			\frac{5\pi^{\frac{3}{2}}\cdot e^{y}}{36\sqrt{3}y^{\frac{3}{2}}}g(y)
		\end{align}
whenever $y\geq 61$, where
\begin{align*}
				f(y)=\left(1-\frac{1}{y^6}\right)\left( D_I(y)-\frac{31}{y^6}\right),
				\qquad g(y)=\left(1+\frac{1}{y^6}\right)\left( D_I(y)+\frac{31}{y^6}\right)
			\end{align*}
and $D_I(y)$ is given by \eqref{defi-I-e}.

For $t\geq 3$, it is routine to verify that
		\begin{align}\label{f0,g0}
			&f(t)>\tilde{f}(t)=1-\frac{3}{8t} - \frac{15}{128t^2} - \frac{105}{1024t^3}-\frac{4725}{32768t^4}-\frac{72765}{262144t^5}-\frac{32}{t^6},\notag\\
			&g(t)<\tilde{g}(t)=1-\frac{3}{8t} - \frac{15}{128t^2} - \frac{105}{1024t^3}-\frac{4725}{32768t^4}-\frac{72765}{262144t^5}+\frac{32}{t^6}.
		\end{align}
Combining \eqref{p6-bound} and \eqref{f0,g0}, we obtain that for $y\geq 61$,
		\begin{align*}
			p_6(n+1)p_6(n)^2-2p_6(n-1)p_6(n+1)^2+p_6(n-1)p_6(n)p_6(n+2)
			\geq\left(\frac{5\pi^{\frac{3}{2}}}{36\sqrt3}\right)^3 F(y),
		\end{align*}
		where
		\begin{align*}
			F(y)=\frac{e^{2y+z}}{{y}^{3}{z}^{\frac{3}{2}}} \tilde{f}(z)\tilde{f}(y)^2
			     -\frac{2e^{x+2z}}{{x}^{\frac{3}{2}}{z}^{3}}\tilde{g}(x)\tilde{g}(z)^2
			     +\frac{e^{x+y+w}}{{x}^{\frac{3}{2}}{y}^{\frac{3}{2}}{w}^{\frac{3}{2}}}\tilde{f}(x)\tilde{f}(y)\tilde{f}(w).
		\end{align*}
		 To establish Theorem \ref{thm-p_kmain}, it is sufficient to show $F(y)>0$.
		Let	
\begin{align}\label{z_6}
\theta_2=\hat{h}_y\left(\frac{5\pi^2}{27}\right).
\end{align}
	We find that the sign of $F(y)$ coincides with that of $F(y)e^{-3\theta_2}$.
	One can verify that for $y\geq 1$,
	\begin{align*}
		2y+z_1-3\theta_2 &= -\frac{25 \left(432 \pi ^4 y^4-160 \pi ^6 y^2+125 \pi ^8\right)}{419904 y^7}<0\\
		x_2+2z_2-3\theta_2 &=- \frac{25 \left(54 \pi ^4 y^2 -5 \pi ^6\right)}{13122 y^5}<0\\
		x_1+y+w_1-3\theta_2 &= -\frac{25 \left(3024 \pi ^4 y^4-1240 \pi ^6 y^2+2125 \pi ^8\right)}{419904 y^7}<0.
	\end{align*}
	Then according to \eqref{eq-Ebound}, we have
	\begin{align*}
		F(y)e^{-3\theta_2}\geq\
		&G_1(2y+z_1-3\theta_2)\frac{\tilde{f}(z)\tilde{f}(y)^2}{{z}^{\frac{3}{2}}{y}^{3}}
		-2G_2(x_2+2z_2-3\theta_2)\frac{\tilde{g}(x)\tilde{g}(z)^2}{{x}^{\frac{3}{2}}{z}^{3}}\\ \nonumber
		&\qquad\qquad\qquad\qquad\qquad\qquad\
		+G_1(x_1+y+w_1-3\theta_2)\frac{\tilde{f}(x)\tilde{f}(y)\tilde{f}(w)}{{x}^{\frac{3}{2}}{y}^{\frac{3}{2}}{w}^{\frac{3}{2}}}.
	\end{align*}
	where $G_1(y)$ and $G_2(y)$ are given by  \eqref{eq-e1-e2}.
Based on \eqref{f0,g0} and \eqref{y bound}, one can verify that  if $y\geq 3$ then 
  {\small
	\begin{align}\label{1,2-bound}
		\tilde{f}(x)&>\lambda_1(y)=1-\frac{3}{8x_1} - \frac{15}{128{x}^2} - \frac{105}{1024{x}^2x_1} - \frac{4725}{32768{x}^4} - \frac{72765}{262144{x}^4x_1} - \frac{32}{{x}^6},\notag\\
		\tilde{f}(z)&>\lambda_2(y)=1-\frac{3}{8z_1} - \frac{15}{128{z}^2} - \frac{105}{1024{z}^2z_1} - \frac{4725}{32768{z}^4} - \frac{72765}{262144{z}^4z_1} - \frac{32}{{z}^6},\notag\\
		\tilde{f}(w)&>\lambda_3(y)=1-\frac{3}{8w_1} - \frac{15}{128{w}^2} - \frac{105}{1024{w}^2w_1} - \frac{4725}{32768{w}^4} - \frac{72765}{262144{w}^4w_1} - \frac{32}{{w}^6},\\
		\tilde{g}(x)&<\lambda_4(y)=1-\frac{3}{8x_2} - \frac{15}{128{x}^2} - \frac{105}{1024{x}^2x_2} - \frac{4725}{32768{x}^4} - \frac{72765}{262144{x}^4x_2} + \frac{32}{{x}^6},\notag\\
		\tilde{g}(z)&<\lambda_5(y)=1-\frac{3}{8z_2} - \frac{15}{128{z}^2} - \frac{105}{1024{z}^2z_2} - \frac{4725}{32768{z}^4} - \frac{72765}{262144{z}^4z_2} + \frac{32}{{z}^6}.\notag
	\end{align}}

	Thus, it is enough to show for $y\geq 61$,
	\begin{align}\label{eq-k-p-E-12}
		&G_1(2y+z_1-3\theta_2)\frac{\lambda_2(y)\tilde{f}(y)^2}{{z}^{\frac{3}{2}}{y}^{3}}
		-2G_2(x_2+2z_2-3\theta_2)\frac{\lambda_4(y)\lambda_5(y)^2}{{x}^{\frac{3}{2}}{z}^{3}}\\ \nonumber
		&\qquad\qquad\qquad\qquad\qquad\qquad\
		+G_1(x_1+y+w_1-3\theta_2)\frac{\lambda_1(y)\tilde{f}(y)\lambda_3(y)}{{x}^{\frac{3}{2}}{y}^{\frac{3}{2}}{w}^{\frac{3}{2}}}>0.
	\end{align}
	Notice that there exist the annoying terms $\sqrt{x}, \sqrt{z}$ and $\sqrt{w}$, we need to do a little change to estimate theses terms.
	Let
	\begin{align}\label{W}
\theta=\sqrt{y^2+\frac{5\pi^2}{27}},\quad W_1=\sqrt{\frac{\theta^9}{y^6 {z}^3}},\quad
		W_2=\sqrt{\frac{\theta^9}{x^3{z}^6}},\quad
		W_3=\sqrt{\frac{\theta^9}{{x}^3 y^3{w}^3}}.
	\end{align}
	Then it can be calculated using Taylor expansion that for $y\geq 5$,
		\begin{align*}
			W_1>W_{11},\quad W_2<W_{22},\quad W_3>W_{31},
		\end{align*}
	where
	\begin{align*}
		W_{11} &= 1 + \frac{25\pi^4 }{324 y^4} - \frac{250 \pi^6 }{6561y^6} + \frac{38125 \pi^8 }{1889568 y^8} - \frac{34375 \pi^{10} }{3188646 y^{10}}, \\
		W_{22} &= 1 + \frac{25\pi^4}{81 y^4} - \frac{250 \pi^6 }{6561 y^6} + \frac{11875 \pi^8 }{118098 y^8}, \\
		W_{31} &= 1 + \frac{175\pi^4}{324 y^4} - \frac{3875 \pi^6 }{13122 y^6}  + \frac{848125 \pi^8 }{1889568 y^8} - \frac{5171875 \pi^{10} }{12754584 y^{10}}.
	\end{align*}
	We find that the left-hand side of \eqref{eq-k-p-E-12} is greater than $\widetilde{F}(y)/\sqrt{\theta^{9}}$, where
	\begin{align*}
		\widetilde{F}(y)=
		&G_1(2y+z_1-3\theta_2)W_{11} \lambda_2(y)\tilde{f}(y)^2
		-2G_2(x_2+2z_2-3\theta_2) W_{22} {\lambda_4(y)\lambda_5(y)^2}\\ \nonumber
		&\qquad\qquad\qquad\qquad\qquad\qquad\
		+G_1(x_1+y+w_1-3\theta_2)W_{31} \lambda_1(y)\lambda_3(y)\tilde{f}(y).
	\end{align*}
By substituting the expressions for $x_1$, $x_2$, $z_1$, $z_2$, $w_1$, $w_2$, and $\theta_2$ into $\widetilde{F}(y)$, we can rewrite $\widetilde{F}(y)$ as
	\begin{align*}
\widetilde{F}(y)=\dfrac{\sum_{k=0}^{104}a_ky^k}{2^{75}3^{38}y^{43}H_1(y)},
	\end{align*}
	where
	{\small\begin{align*}
		H_1(y)=&2^{75}3^{38}y^{43}
		\left(9 y^2-5 \pi ^2\right)^3 \left(9 y^2+5 \pi ^2\right)^6
		\left(9 y^2+10 \pi ^2\right)^3\\
		&\times\left(419904 y^8+116640 \pi ^2 y^6-16200 \pi ^4 y^4+4500 \pi ^6 y^2-3125 \pi ^8\right)\\
		&\times\left(419904 y^8-116640 \pi ^2 y^6-16200 \pi ^4 y^4-4500 \pi ^6 y^2-3125 \pi ^8\right)\\
		&\times
		\left(26244 y^8+14580 \pi ^2 y^6-4050 \pi ^4 y^4+2250 \pi ^6 y^2-3125 \pi ^8\right)\\
		&\times{\left(11664 y^6+3240 \pi ^2 y^4-450 \pi ^4 y^2+125 \pi ^6\right)^2 }
		\left(11664 y^6-3240 \pi ^2 y^4-450 \pi ^4 y^2-125 \pi ^6\right),
	\end{align*} }
	and $a_k$ are the known numbers. Specially, we give the values of $a_{104}$ and $a_{103}$ below:
	{\begin{align*}
		a_{104}=2^{98}3^{99}5^{3}\pi^6,\qquad a_{103}=2^{95}3^{96}(2^{2}5^{4}\pi^8-3^{3}5^{3}41\pi^6-2^{13}3^{9}).
	\end{align*}}
One can check that $H_1(y)$ is positive for $y\ge 3$.
It remains to show that
	\begin{align*}
		H(y)=\sum_{k=0}^{104}a_ky^k>0.
	\end{align*}
It can be computed by mathematical software that the largest real zero of $H(y)$ is less than $12$.
Thus, for $y\geq 61$, i.e, $n\geq 679$,
 we have $F(y)>0$ along with the fact that $H(12)>0$. Additionally, for $60\leq n\leq 678$ one can directly verify \eqref{p_kmain}. This completes the proof.
\end{proof}

For the log-concavity of the $k$-regular partition function, we have the following result.

\begin{thm}\label{lc-k-regular-pf}
For $2\leq k\leq 9$, let $N_k$ be given as in Theorem \ref{thm-p_kmain}. Then the sequence $\{p_k(n)\}_{n\ge N_k}$ is log-concave.
\end{thm}

\begin{proof}
The log-concavity of $\{p_k(n)\}_{n\ge 58}$ when $k=2,3,4$ or $5$ has been proved by Dong and Ji in \cite[Theorem 1.4]{Dong-Ji-2024}. For each $6\le k\le 9$, one can give a proof the log-concavity of $\{p_k(n)\}_{n\ge 36}$,  exactly like that of Theorem \ref{thm-p_kmain}.
Again we take $k=6$ to illustrate the idea,
After some analysis, we obtain that for $y\geq 61$,
\begin{align*}
p_6(n)^2-p_6(n-1)p_6(n+1)\geq \left(\frac{5\pi^{\frac{3}{2}}}{36\sqrt3}\right)^2 e^{2y}J(y),
\end{align*}
where
\begin{align*}
	J(y)=\frac{1}{{y}^{3}}\tilde{f}(y)^2
	-\frac{e^{x_2+z_2-2y}}{{x}^{\frac{3}{2}}{z}^{\frac{3}{2}}}\tilde{g}(x)\tilde{g}(z),
\end{align*}
the symbols $x,y,z,x_2,z_2$ are defined as in \eqref{eq-y_k} and \eqref{+--bound}, and the functions $\tilde{f}(t),\,\tilde{g}(t)$ are given by \eqref{f0,g0}. It can be verified that $x_2+z_2-2y<0$ for $y\ge 1$. According to \eqref{eq-Ebound} and the above bounds of $\tilde{g}(x)$ and $\tilde{g}(z)$ (immediately before \eqref{eq-k-p-E-12}), we have
\begin{align}\label{J_1}
	J(y)\ge \frac{1}{{y}^{3}}\tilde{f}(y)^2
	-\frac{G_2(x_2+z_2-2y)}{{x}^{\frac{3}{2}}{z}^{\frac{3}{2}}}\lambda_4(y)\lambda_5(y),
\end{align}
which $\lambda_4(y)$ and $\lambda_5(y)$ are defined in \eqref{1,2-bound}.
Let
$$V=\sqrt{\frac{y^6}{{x}^{3}{z}^{3}}}\qquad\text{and}\qquad V_2=1+\frac{25\pi^4}{108y^4}+\frac{4375\pi^8}{69984y^8}+\frac{240625\pi^{12}}{12754584 y^{12}}.$$
Using Taylor expansion, it can be checked that $V<V_2$ for $y\ge 4$. Thus the right-hand side of \eqref{J_1} is greater than $\widetilde{J}(y)/y^3$, where
\begin{align*}
\widetilde{J}(y)=\tilde{f}(y)^2
-G_2(x_2+z_2-2y) V_2 {\lambda_4(y)\lambda_5(y)}.
\end{align*}
After some simplification, we see that
\small\begin{align*}
	\widetilde{J}(y)=\dfrac{\sum_{l=0}^{39}a_ly^l}{J_1(y)},
\end{align*}
where {\small
	\begin{align*}
	J_1(y)=&2^{40}3^{21}(9y^2-5\pi^2)^3(9y^2+5\pi^2)^3(11664y^6-3240\pi^2y^4-450\pi^4y^2-125\pi^6)\\
	&\times(11664y^6+3240\pi^2y^4-450\pi^4y^2+125\pi^6),
	\end{align*}} and
$a_l$ are the known number, and the values of $a_{39}$ and $a_{38}$ are given below
\begin{align*}
	a_{39}=2^{46}3^{41}5^2\pi^4,\quad a_{38}=-2^{44}3^{42}5^3\pi^4.
\end{align*}
One can check that $J_1(y)$ is positive for $y\ge 3$. Then it is sufficient to show the denominator of $\widetilde{J}(y)>0$.
It can be computed by mathematical software that the largest real zero of $\widetilde{J}(y)$ is less than $5$.
Thus, for $y\geq 61$, i.e, $n\geq 679$, we have $\widetilde{J}(y)>0$. Additionally, for $36\leq n\leq 678$ one can directly verify that $p_6(n)^2-p_6(n-1)p_6(n+1)\geq 0$. This completes the proof.
\end{proof}

Based on Theorems \eqref{thm-p_kmain} and \eqref{lc-k-regular-pf}, we immediately obtain the following result.

\begin{cor}
For $2\leq k\leq 9$, let $N_k$ be given as in Theorem \ref{thm-p_kmain}. Then the sequence $\{p_k(n)\}_{n\ge N_k}$ satisfies the Briggs inequality.
\end{cor}

\subsection{$k$-regular overpartition}

For the $k$-regular overpartition function $\overline{p}(n)$, we obtain the following result.

\begin{thm}\label{thm-opmain}
	For $2\leq k\leq 9$ and $n\geq \overline{N}_k$, we have
	\begin{align}\label{opk-main}
		\overline{p}_k(n+1)\overline{p}_k(n)^2-2\overline{p}_k(n-1)\overline{p}_k(n+1)^2+\overline{p}_k(n-1)\overline{p}_k(n)\overline{p}_k(n+2)>0,
	\end{align}
	where
	\begin{align*}
		&\overline{N}_2=75,\quad \overline{N}_3=17,\quad \overline{N}_4=33,\quad \overline{N}_5=30,\quad\\
		& \overline{N}_6=10,\quad \overline{N}_7=24,\quad \overline{N}_8=27,\quad \overline{N}_9=10.
	\end{align*}
\end{thm}

To prove Theorem \ref{thm-opmain}, we need the following 
upper and lower bounds of $\p_k(n)$ for $2\le k\le 9$ given by
Peng, Zhang and Zhong  \cite{Peng-Zhang-Zhong-2023}.

\begin{thm}[\cite{Peng-Zhang-Zhong-2023}, Corollary 3.4]\label{thm-op}
	Let $I_1(s)$ denote the first modified Bessel function of the first kind and
	$$\overline{\mu}_k(n)=\pi\sqrt{\left(1-\frac{1}{k}\right)n}.$$
	For $2\leq k\leq  9$ and $n\geq \overline{n}_k$, we have
	\begin{align*}
		{\M_k(n)}\left(1-\frac{1}{{\overline{\mu}_k(n)}^6}\right)\leq \p_k(n)\leq {\M_k(n)}\left(1+\frac{1}{{\overline{\mu}_k(n)}^6}\right),
	\end{align*}
where $\M_k(n)=C_k(n)I_1(\bmu_k(n))$, and the values of $\overline{n}_k$ and $C_k(n)$ are given in Table \ref{opk(n)-bounds}. 
	\begin{table}[h]
		\centering	
		\begin{tabular}{c|c|c|c|c|c|c|c|c}
			$k$ & $2$ & $3$ &$4$ & $5$ & $6$&$7$&$8$&$9$
			\\[3pt]   \hline \rule{0pt}{25pt}
			$C_k(n)$&$\frac{\pi^2}{\sqrt{8}\bmu_2(n)}$&$\frac{2\sqrt{3}\pi^2}{9\bmu_3(n)}$&$\frac{3\pi^2}{4\bmu_4(n)}$&$\frac{8\sqrt{5}\pi^2}{25\bmu_5(n)}$&$\frac{5\sqrt{6}\pi^2}{18\bmu_6(n)}$&$\frac{18\sqrt{7}\pi^2}{49\bmu_7(n)}$&$\frac{7\sqrt{2}\pi^2}{8\bmu_8(n)}$&$\frac{8\pi^2}{9\bmu_9(n)}$
			\\[3pt]   \hline \rule{0pt}{23pt}
			$\overline{n}_k$&$375$&$365$&$250$&$427$&$2055$&$1230$&$1927$&$8187$
			\\[3pt]
		\end{tabular}
		\caption{Values of $\overline{n}_k$ and $C_k(n)$ for $2\leq k\leq9$.}
		\label{opk(n)-bounds}
	\end{table}
\end{thm}

Now we are in the position to prove Theorem \ref{thm-opmain}. 

\begin{proof}[Proof of Theorem \ref{thm-opmain}]
		We shall take $k=6$ to illustrate our proof, and the proofs for other values of $k$ can be given in the same manner. Set
		\begin{align}\label{op-x}
			x=\bmu_6(n-1),\quad y=\bmu_6(n),\quad
			z=\bmu_6(n+1),\quad
			w=\bmu_6(n+2).
		\end{align}
		Then
		\begin{align*}	
			x=\sqrt{y^2-\frac{5\pi^2}{6}},\quad  z=\sqrt{y^2+\frac{5\pi^2}{6}},\quad w=\sqrt{y^2+\frac{5\pi^2}{3}}.
		\end{align*}
		By Theorem \ref{thm-op}, we find that for $n\ge 2055$, i.e., $y\ge 130$,
		\begin{align*}
			{\frac{5\sqrt{6}\pi^2}{18y}}I_1(y)\left(1-\frac{1}{{y}^6}\right)< \p_6(n)< {\frac{5\sqrt{6}\pi^2}{18y}}I_1(y)\left(1+\frac{1}{{y}^6}\right).
		\end{align*}
		Applying Lemma \ref{Ji-Dong}, we have
		\begin{align}\label{op6-bound}
			\frac{5\pi^{\frac{3}{2}}\cdot e^{y}}{6\sqrt{3}y^{\frac{3}{2}}}f(y)
			<\p_6(n)<
			\frac{5\pi^{\frac{3}{2}}\cdot e^{y}}{6\sqrt{3}y^{\frac{3}{2}}}g(y)
		\end{align}
		whenever $y\geq 130$, where the functions $f(y),g(y)$ are given by \eqref{f0,g0}. Using the symbols given by \eqref{+--bound} and \eqref{z_6}, we set
		\begin{align}\label{xzw-bound}
			x_1=\check{h}_{y}\left(-\frac{5\pi^2}{6}\right),\quad
			z_1=\check{h}_{y}\left(\frac{5\pi^2}{6}\right),\quad
			w_1=\check{h}_{y}\left(\frac{5\pi^2}{3}\right),\notag\\
			x_2=\hat{h}_{y}\left(-\frac{5\pi^2}{6}\right),\quad
			z_2=\hat{h}_{y}\left(\frac{5\pi^2}{6}\right),\quad
			w_2=\hat{h}_{y}\left(\frac{5\pi^2}{3}\right),
		\end{align}
		and $\theta_2=\hat{h}_{y}\left(\frac{5\pi^2}{18}\right)$.
		
		Similar to the proof of Lemma \ref{lem-x-estimation}, we can show that if $y\geq 3$ then
		\begin{align*}
			x_1<x<x_2,\quad
			z_1<z<z_2,\quad
			w_1<w<w_2.
		\end{align*}
		Recall the expressions of $\tilde{f}(t)$ and $\tilde{g}(t)$ in \eqref{f0,g0}. Combining \eqref{op6-bound} and \eqref{f0,g0}, we obtain that for $y\geq 130$,
		\begin{small}
		\begin{align*}
			\p_6(n+1)\p_6(n)^2-2 \p_6(n-1)\p_6(n+1)^2+\p_6(n-1)\p_6(n)\p_6(n+2)
			\ge\left(\frac{5\pi^{\frac{3}{2}}}{6\sqrt3}\right)^3 F(y),
		\end{align*}
	\end{small}
		where
		\begin{align*}
			F(y)=\frac{e^{2y+z}}{{y}^{3}{z}^{\frac{3}{2}}} \tilde{f}(z)\tilde{f}(y)^2
			-\frac{2e^{x+2z}}{{x}^{\frac{3}{2}}{z}^{3}}\tilde{g}(x)\tilde{g}(z)^2
			+\frac{e^{x+y+w}}{{x}^{\frac{3}{2}}{y}^{\frac{3}{2}}{w}^{\frac{3}{2}}}\tilde{f}(x)\tilde{f}(y)\tilde{f}(w).
		\end{align*}
		A little computation shows that $2y+z_1-3\theta_2<0,\, x_2+2z_2-3\theta_2<0$ and $x_1+y+w_1-3\theta_2 <0$ for $y\geq 2$. By \eqref{eq-Ebound}, we have
		\begin{align*}
			F(y)e^{-3\theta_2}\geq
			&G_1(2y+z_1-3\theta_2)\frac{\tilde{f}(z)\tilde{f}(y)^2}{{z}^{\frac{3}{2}}{y}^{3}}
			-2G_2(x_2+2z_2-3\theta_2)\frac{\tilde{g}(x)\tilde{g}(z)^2}{{x}^{\frac{3}{2}}{z}^{3}}\\ \nonumber
			&\qquad\qquad\qquad\qquad\qquad\qquad\
			+G_1(x_1+y+w_1-3\theta_2)\frac{\tilde{f}(x)\tilde{f}(y)\tilde{f}(w)}{{x}^{\frac{3}{2}}{y}^{\frac{3}{2}}{w}^{\frac{3}{2}}}.
		\end{align*}
		Using the same arguments as in the proof of Theorem \ref{thm-p_kmain}, it is sufficient to show for $y\geq 130$,
		\begin{align}\label{opE}
			&G_1(2y+z_1-3\theta_2)\frac{\lambda_2(y)\tilde{f}(y)^2}{{z}^{\frac{3}{2}}{y}^{3}}
			-2G_2(x_2+2z_2-3\theta_2)\frac{\lambda_4(y)\lambda_5(y)^2}{{x}^{\frac{3}{2}}{z}^{3}}\\ \nonumber
			&\qquad\qquad\qquad\qquad\qquad\qquad\
			+G_1(x_1+y+w_1-3\theta_2)\frac{\lambda_1(y)\tilde{f}(y)\lambda_3(y)}{{x}^{\frac{3}{2}}{y}^{\frac{3}{2}}{w}^{\frac{3}{2}}}>0.
		\end{align}
		Using the symbols given by \eqref{W}, we set
		\begin{align*}
			\theta=\sqrt{y^2+\frac{5\pi^2}{18}},\quad W_1=\sqrt{\frac{\theta^9}{y^6 {z}^3}},\quad
			W_2=\sqrt{\frac{\theta^9}{x^3{z}^6}},\quad
			W_3=\sqrt{\frac{\theta^9}{{x}^3 y^3{w}^3}}.
		\end{align*}
		and
		\begin{align*}
			W_{11} &= 1 + \frac{25\pi^4 }{144y^4} - \frac{125 \pi^6 }{972y^6} + \frac{38125 \pi^8 }{373248 y^8} - \frac{34375 \pi^{10} }{419904 y^{10}}, \\
			W_{22} &= 1 + \frac{25\pi^4}{36 y^4} - \frac{125 \pi^6 }{972 y^6} + \frac{11875 \pi^8 }{23328 y^8}, \\
			W_{31} &= 1 + \frac{175\pi^4}{144 y^4} - \frac{3875 \pi^6 }{3888 y^6}  + \frac{848125 \pi^8 }{373248 y^8} - \frac{5171875 \pi^{10} }{1679616 y^{10}}.	
		\end{align*}
		One can check that for $y\ge 5$,
		$
		W_1>W_{11}, W_2<W_{22}, W_3>W_{31}.
		$
		We find that the left-hand side of \eqref{opE} is greater than $\widetilde{F}(y)/\sqrt{\theta^{9}}$, where
		\begin{align*}
			\widetilde{F}(y)=
			&G_1(2y+z_1-3\theta_2)W_{11} \lambda_2(y)\tilde{f}(y)^2
			-2G_2(x_2+2z_2-3\theta_2) W_{22} {\lambda_4(y)\lambda_5(y)^2}\\ \nonumber
			&\qquad\qquad\qquad\qquad\qquad\qquad\
			+G_1(x_1+y+w_1-3\theta_2)W_{31} \lambda_1(y)\lambda_3(y)\tilde{f}(y).
		\end{align*}
	Note that for $1\le i\le 5$, the function $\lambda_i(y)$ is defined in a similar way to \eqref{1,2-bound} by substituting the variables in \eqref{1,2-bound} with the variables in \eqref{op-x} and  \eqref{xzw-bound}.
		After some simplification, we see that
		\small\begin{align*}
			\widetilde{F}(y)=\dfrac{\sum_{l=0}^{104}a_ly^l}{H_1(y)},
		\end{align*}
		where
	{\small	\begin{align*}
			H_1(y)=&2^{89}3^{24}y^{43} \left(5 \pi^2 + 3 y^2\right)^3 \left(-5 \pi^2 + 6 y^2\right)^3 \left(5 \pi^2 + 6 y^2\right)^6\\
&\times\left(-125 \pi^6 - 300 \pi^4 y^2 - 1440 \pi^2 y^4 + 3456 y^6\right)\left(125 \pi^6 - 300 \pi^4 y^2 + 1440 \pi^2 y^4 + 3456 y^6\right)^2\\
&\times \left(-3125 \pi^8 + 1500 \pi^6 y^2 - 1800 \pi^4 y^4 + 4320 \pi^2 y^6 + 5184 y^8\right)
			 \\
			&\times \left(-3125 \pi^8 - 3000 \pi^6 y^2 - 7200 \pi^4 y^4 - 34560 \pi^2 y^6 + 82944 y^8\right)\\
			&\times\left(-3125 \pi^8 + 3000 \pi^6 y^2 - 7200 \pi^4 y^4 + 34560 \pi^2 y^6 + 82944 y^8\right),
	    \end{align*}}
		and $a_l$ are the known numbers, and the values of $a_{104}$ and $a_{103}$ are given below
		\begin{align*}
			a_{104}=2^{139}3^{55}5^3\pi^6,\quad a_{103}=2^{136}3^{53}(2\cdot 5^4\pi^8-3^25^341\pi^6-2^{16}3^5).
		\end{align*}
		One can check that $H_1(y)$ is positive for $y\ge 3$. Then it is sufficient to show the numerator of $\widetilde{F}(y)>0$.
		It can be computed by mathematical software that the largest real zero of $\widetilde{F}(y)$ is less than $9$.
		Thus, for $y\geq 130$, i.e, $n\geq 2055$, we have $\widetilde{F}(y)>0$. Additionally, for $10\leq n\leq 2054$ one can directly verify that inequality \eqref{opk-main}. This completes the proof.
	\end{proof}

For $2\le k\le 9$, Peng, Zhang and Zhong \cite{Peng-Zhang-Zhong-2023} proved that the log-concavity of  $\{p_k(n)\}_{n\ge 21}$. Based on their result and Theorem \ref{thm-opmain}, we have the following result.

\begin{thm}
For $2\leq k\leq 9$, let $\overline{N}_k$ be given as in Theorem \ref{thm-opmain}. Then the sequence $\{\p_k(n)\}_{n\ge \overline{N}_k}$ satisfies Briggs inequality.	
\end{thm}

Note that the conditions of Chern's formula in \cite{Chern-2019} are not suitable for $k\ge 25$ when considering $p_k(n)$ or for $k\ge 10$ when considering $\p_k(n)$, thus we cannot prove the inequalities \eqref{p_kmain} or \eqref{opk-main} for any $k$. But we check that for $2\le k\le 100$, \eqref{p_kmain} stands when $220\le n\le 1000$, while \eqref{opk-main} stands when $75\le n\le 1000$. Hence, we propose the following conjecture.
\begin{conj}\label{p_k_cj}
	For any positive integer $k\geq 2$, there exist $N(k)$ such that the inequalities \eqref{p_kmain} and \eqref{opk-main} hold for $n\geq N(k)$.
\end{conj}
\vspace{10pt}
\section{Appendix: proof of Theorem \ref{thm-p_k2}}
In this appendix, we will follow Dong and Ji \cite{Dong-Ji-2024} to establish upper and lower bounds of $p_k(n)$ for $3\leq k\leq 9$ by using Chern's formula \cite{Chern-2019} for $\eta$-quotients.

Firstly, let us review Chern's theorem.
We adopt the notations in \cite{Dong-Ji-2024}.
Let $\mathbf{m}=(m_1,\ldots,m_R)$ be a sequence of $R$ distinct positive integers and $\mathbf{\delta}=(\delta_1,\ldots,\delta_R)$ be a sequence of $R$ non-zero integers.
Assuming that $h$ and $j$ are positive integers with ${\rm gcd}(h, j) =1$,  set
\[C_1=-\frac{1}{2}\sum_{r=1}^R\delta_r,\qquad C_2=\sum_{r=1}^Rm_r\delta_r,\]
\[C_3(l)=-\sum_{r=1}^R\frac{\delta_r\gcd^2(m_r,l)}{m_r},\qquad C_4(l)=\prod_{r=1}^R\left(\frac{m_r}{\gcd(m_r,l)}\right)^{-\frac{\delta_r}{2}},\]
\begin{equation*}
	\hat{A}_{l}(n)=\sum_{0\leq h< l\atop {\rm gcd}(h,l)=1}\exp\left(-\frac{2\pi nh i }{l}-\pi i\sum_{r=1}^R\delta_rs\left(\frac{m_rh}{\gcd(m_r,l)},\frac{l}{\gcd(m_r,l )}\right)\right),
\end{equation*}
where  $s(h,j)$ is the  Dedekind sum.
Take $L=\mathrm {lcm} (m_1,\ldots,m_R)$, the least common multiple of $m_1,\ldots,m_R$. We divide the set  $\{1,2\cdots,L\}$ into the following two disjoint subsets:
\begin{align*}
	\mathcal{L}_{>0}:=\{1\le l\le L\ |\ C_3(l)>0\},\qquad
	\mathcal{L}_{\le0}:=\{1\le l\le L\ | \ C_3(l)\le0\}.
\end{align*}
Define
\begin{align*}
	\mathbb{E}_{C_1}(s):=
	\left\{\begin{array}{ll}
		1,~~~~~~~~~~~~~~~~~~&C_1=0,\\[5pt]
		2\sqrt{s},~~~~~~~~~~~&C_1=-\frac{1}{2},\\[5pt]
		s\log(s+1),~~~~~~~&C_1=-1,\\[5pt]
		s^{-2C_1-1}\zeta(-C_1),~~~~~~~&\mbox{otherwise},\\[5pt]
	\end{array}\right.
\end{align*}
where $\zeta(\cdot)$ is Riemann zeta-function.
Now, we are able to give Chern's theorem.
\begin{thm}[\cite{Chern-2019}, Theorem 1.1]\label{lem-C} 	
	Let
$$G(q)=\sum_{n\geq 0} g(n)q^n=\prod_{r=1}^R(q^{m_r};q^{m_r})_{\infty}^{\delta_r}.$$ If  $C_1\leq0$ and the inequality
	\begin{equation}\label{eq-min}
		\min_{1\le r\le R}\left(\frac{\gcd^2(m_r,l)}{m_r}\right)\ge \frac{C_3(l)}{24}
	\end{equation}
	holds for all  $1\le l \le L$, then for positive integers $N$ and  $n>-\frac{c_2}{24}$, we have
	\begin{small}
	\begin{align}\label{eq-main-1}
		g(n)=E(n)+&\sum_{l\in \mathcal{L}_{>0}}2\pi C_4(l)\left(\frac{24n+C_2}{C_3(l)}\right)^{-\frac{C_1+1}{2}}
		\sum_{\substack{1\le t\le N\\t\equiv_L l }} {I_{-C_1-1}\left(\frac{\pi}{6t}\sqrt{C_3(l)(24n+C_2)}\right)} \frac{\hat{A}_{t}(n)}{t},
	\end{align}
	\end{small}
	where
	\begin{align*}
		|E(n)|
		&\le \frac{2^{-C_1}\pi^{-1}N^{-C_1+2}}{n+\frac{C_2}{24}}
		\exp\left(\frac{2\pi}{N^2} \left(n+\frac{C_2}{24}\right)\right) \sum_{l\in \mathcal{L}_{>0}} \exp\left(\frac{C_3(l)\pi}{3}\right)\\[5pt]
		&\quad + 2 \exp\left(\frac{2\pi}{N^2} \left(n+\frac{C_2}{24}\right)\right) \mathbb{E}_{C_1}(N) \times \Bigg(-\sum_{l\in \mathcal{L}_{>0}}C_4(l)\exp\left(\frac{\pi C_3(l)}{24}\right)\\[5pt]
		&\quad + \sum_{1\le l\le L}C_4(l)
		\exp
		\bigg(\frac{\pi C_3(l)}{24}+\sum_{r=1}^R\frac{|\delta_r|\exp\left(-\pi\gcd^2(m_r,l)/m_r\right)}{\left(1-\exp\left(-\pi\gcd^2(m_r,l)/m_r\right)\right)^2}\bigg)\Bigg).
	\end{align*}
	and $I_{\nu}(s)$ is the $\nu$-th modified Bessel function of the first kind.
\end{thm}

Recall that, for $k\geq 2$, the generating function for the sequence $\{p_k(n)\}_n$ is as follows:
	\begin{align*}
		\sum_{n\geq 0}p_k(n)q^n=\prod_{n=1}^{\infty}\dfrac{1-q^{kn}}{1-q^n}=\frac{(q^k;q^k)_\infty}{(q;q)_\infty},
	\end{align*}
	where $(a;q)_\infty=\prod_{j=0}^\infty (1-aq^j)$.
In order to use Theorem \ref{lem-C} to get bounds of $p_k(n)$, we need to $g(n)=p_k(n)$ in Chern's theorem,
in which case we have $\mathbf{m}=(1,k)$, $\delta=(-1,1)$,  $R=2$, $L=k$ and
\begin{align}\label{eq-delta1234}
	C_1=0,\qquad C_2=k-1,\qquad
	C_3(l)=1-\frac{\gcd(k,l)^2}{k},\qquad C_4(l)=\sqrt{\frac{{\gcd(k,l)}}{k}}.
\end{align}
We also need to check the conditions of theorem \ref{lem-C}.
The condition $C_1\leq 0$ naturally holds for any $k\geq 2$, while inequality \eqref{eq-min}  is satisfied only for  $2\leq k\leq 24$.
Consequently, it is feasible to apply Chern's theorem to asymptotically estimate $p_k(n)$ for  $2\leq k\leq 24$.
We have the following result.

\begin{prop}\label{prop-k-n-24}
For any $k\geq 2$, let $c_2(k)=k-1, c_3(k,l)=1-\frac{\gcd(k,l)^2}{k},c_4(k,l)=\sqrt{\frac{{\gcd(k,l)}}{k}}$ and
	\begin{align*}
		\mu_k(n)=\frac{\pi}{6}\sqrt{c_3(k,1)(24n+c_2(k))}.
	\end{align*}
	Then, for $2\leq k\leq 24$ and any positive integers $n, N$, we have
	\begin{align}\label{eq-pkn-formula}
p_k(n)=\frac{(k-1)\pi^2}{3k\sqrt{k}\mu_k(n)}I_1(\mu_k(n))+R_k(n),
	\end{align}
	where $R_k(n)= E_k(n) + B_k(n)$ and
	\begin{footnotesize}
	\begin{align}\label{ieq-Ek}
		|E_k(n)|
		&\le\frac{24N^{2}}{\pi(24n+c_2(k))}\exp\left(\frac{\pi\left(24n+c_2(k)\right)}{12N^2} \right)\sum_{l\in \mathcal{L}_{>0}}\exp\left(\frac{c_3(k,l)\pi}{3}\right) +2\exp\left(\frac{\pi\left(24n+c_2(k)\right)}{12N^2} \right) \nonumber \\[5pt]
		&\quad
		\times\Bigg(\sum_{1\le l\le k}c_4(k,l)
		\exp\bigg(\frac{\pi c_3(k,l)}{24} + \frac{\exp(-\pi)}{(1-\exp(-\pi))^2}
		+\frac{\exp\left(-\pi\gcd^2(k,l)/k\right)}{\left(1-\exp\left(-\pi\gcd^2(k,l)/k\right)\right)^2}\bigg) \nonumber \\[5pt]
		&\qquad\qquad \qquad - \sum_{l\in \mathcal{L}_{>0}}c_4(k,l)\exp\left(\frac{\pi c_3(k,l)}{24}\right)\Bigg).
	\end{align}
	\end{footnotesize}
	and
	\begin{align}\label{eq-bkn}
		|B_k(n)|\leq	\frac{\pi^2c_4(k,l')c_3(k,1)}{3\mu_k(n)}
		\sum_{l\in \mathcal{L}_{>0}}
		\sum_{\substack{{2\le t\le N}\\t\equiv_k l }}I_1\left(\frac{\mu_k(n)}{t}\right),
	\end{align}
	with
	$c_4(k,l')=\max\{ c_4(k,l) \ |\ l\in \mathcal{L}_{>0}\}$.
\end{prop}
\begin{proof}
	By substituting the values of \eqref{eq-delta1234} into Theorem \ref{lem-C}, we get
	\begin{small}
	\begin{align}\label{eq-p-k-n}
		p_k(n)=E_k(n)+\sum_{l\in \mathcal{L}_{>0}}2\pi \frac{c_4(k,l)\sqrt{c_3(k,l)}}{\sqrt{24n+c_2(k)}}\times
		\sum_{\substack{1\le t\le N\\t\equiv_k l }}I_{1}\left(\frac{\pi}{6t}\sqrt{c_3(k,l)(24n+c_2(k))}\right)\frac{\hat{A}_{t}(n)}{t},
	\end{align}
	\end{small}
	where $E_k(n)$ plays the role of $E(n)$ of \eqref{eq-main-1} and satisfies \eqref{ieq-Ek}.
	We find that for any $2\leq k\leq 24$, $1\in \mathcal{L}_{>0}$. Moreover, $1\equiv_k 1$ for any $k$, and hence $t=1$ will appear in the second sumnation of \eqref{eq-p-k-n}.
	Based on the fact that $\hat{A}_{1}(n)=1$, \eqref{eq-p-k-n} can be written as
	\begin{align*}
		p_k(n)
		=\frac{\pi^2c_4(k,1)c_3(k,1)}{3\mu_k(n)}I_1(\mu_k(n)) + E_k(n)+B_k(n),
	\end{align*}
	where
	\begin{align*}
		B_k(n)
		=&\frac{\pi^2c_4(1)c_3(1)}{3\mu_k(n)}\sum_{\substack{2\le t\le N\\t\equiv_k 1 }}I_1\left(\frac{\mu_k(n)}{t}\right)\frac{\hat{A}_{t}(n)}{t}\\
		&+\sum_{l\in \mathcal{L}_{>0}\setminus\{1\}}\frac{\pi^2c_4(k,l)\sqrt{c_3(k,l)}\sqrt{c_3(k,1)}}{3\mu_k(n)}\times
		\sum_{\substack{{2\le t\le N}\\t\equiv_k l }}I_1\left(\frac{\sqrt{c_3(k,l)}\mu_k(n)}{\sqrt{c_3(k,1)}t}\right)\frac{\hat{A}_{t}(n)}{t}.
	\end{align*}
	Noting that $c_3(k,l)\leq c_3(k,1)$ for $l\geq 1$, $I_1(s)$ is increasing for $s>0$ and $|\hat{A}_t(n)|\leq t$ for $t\geq 1$, we get
	\begin{align*}
		|B_k(n)|\leq \frac{\pi^2c_4(1)c_3(1)}{3\mu_k(n)}\sum_{\substack{2\le t\le N\\t\equiv_k 1 }}I_1\left(\frac{\mu_k(n)}{t}\right)
		+\sum_{l\in \mathcal{L}_{>0}\setminus\{1\}}\frac{\pi^2c_4(l) {c_3(1)}}{3\mu_k(n)}\times
		\sum_{\substack{{2\le t\le N}\\t\equiv_k l }}I_1\left(\frac{\mu_k(n)}{t}\right).
	\end{align*}
We immediately obtain the desired result.
\end{proof}

One can see that the above bound of $\hat{R_k}(n)$ in \eqref{eq-pkn-formula}
is complicated, and it is not sufficient for our purpose. By using the following upper bound on the first modified Bessel function of the first kind $I_1(s)$:
	\begin{align}\label{eq-lem-B-K}
		I_1(s)\le \sqrt{\frac{2}{\pi s}}e^s,
	\end{align}
due to Bringmann, Kane, Rolen and Trippin \cite{Bringmann-Kane-Rolen-Tripp-2021}, we are able to give a simpler bound of $\hat{R_k}(n)$. We have the following result.

\begin{thm}\label{thm-p_k1}
For $2\leq k\leq 9$, let $\mu_k(n)$ and $R_k(n)$ be given as in Proposition \eqref{prop-k-n-24}, and let $\hat{n}_k$ and $\hat{R_k}(n)$ be given as in Table \ref{pk(n)-bounds}.
	\begin{table}[h]
		\centering
		\begin{tabular}{c|c|c||c|c|c}
			$k$ & $\hat{n}_k$ & $\hat{R}_k(n)$ &$k$ & $\hat{n}_k$ & $\hat{R}_k(n)$
			\\[3pt]   \hline \rule{0pt}{23pt}
			$2$ & $60$
			& {\large$\frac{\pi^{\frac{3}{2}}}{3\sqrt 2\sqrt{\mu_2(n)}}$}
			$\exp\left(\frac{{\mu_2(n)}}{2}\right)$
			&$6$ & $46$
			& {\large$\frac{20\pi^{\frac{3}{2}}}{27\sqrt{3\mu_6(n)}}$} $\exp\left(\frac{\mu_6(n)}{2}\right)$
			\\[9pt] \hline \rule{0pt}{23pt}
			$3$ & $40$
			& {\large$\frac{16\pi^{\frac{3}{2}}}{27\sqrt{3\mu_3(n)}}$}
			$\exp\left(\frac{\mu_3(n)}{2}\right)$
			&$7$ & $75$
			& {\large$\frac{48\pi^{\frac{3}{2}}}{49\sqrt{7\mu_7(n)}}$}
			$\exp\left(\frac{\mu_7(n)}{2}\right)$
			\\[9pt] \hline \rule{0pt}{23pt}
			$4$ & $42$
			& {\large$\frac{\pi^{\frac{3}{2}}}{4\sqrt{\mu_4(n)}}$}
			$\exp\left(\frac{\mu_4(n)}{2}\right)$
			&$8$ & $87$
			& {\large$\frac{7\pi^{\frac{3}{2}}}{16\sqrt{\mu_8(n)}}$}
			$\exp\left(\frac{\mu_8(n)}{2}\right)$
			\\[9pt] \hline \rule{0pt}{23pt}
			$5$ & $47$
			& {\large$\frac{64\pi^{\frac{3}{2}}}{75\sqrt{5\mu_5(n)}}$}
			$\exp\left(\frac{\mu_5(n)}{2}\right)$
			&$9$ & $130$
			& {\large$\frac{64\pi^{\frac{3}{2}}}{243\sqrt{\mu_9(n)}}$}
			$\exp\left(\frac{\mu_9(n)}{2}\right)$
			\\[3pt]
		\end{tabular}
		\caption{Values of $\hat{n}_k$ and $\hat{R}_k(n)$ for $2\leq k\leq9$.}
		\label{pk(n)-bounds}
	\end{table}
If $n\geq \hat{n}_k$, then $|R_k(n)|\leq \hat{R_k}(n)$. \end{thm}

\begin{proof}
We prove the theorem for $k=6$ and omit the details for other values of $k$. 
By Proposition \ref{prop-k-n-24}, it is enough to show that for $n\ge 46$,
	\begin{align*}
		|E_6(n)| + |B_6(n)| \leq \frac{20\pi^{\frac{3}{2}}}{27\sqrt{3\mu_6(n)}},
	\end{align*}
where $E_6(n)$ satisifies \eqref{ieq-Ek} and $B_6(n)$ satisfies \eqref{eq-bkn}.

To further estimate $E_6(n)$ and $B_6(n)$, 
we first determine the values of $c_3(6,l)$ and $c_4(6,l)$, which are listed in Table \ref{tab:1}.
	\begin{table}[htbp]
		\centering
		\begin{tabular}{ccccccc}
			\hline\noalign{\smallskip}
			$l$ & 1 & 2 & 3 & 4 & 5 & 6\\
			\noalign{\smallskip}\hline\noalign{\smallskip}
			$c_3(6,l)$  & $ \frac{5}{6}$ &  $\frac{1}{3}$ &  $-\frac{1}{2}$ &  $\frac{1}{3}$ & $\frac{5}{6}$ & $-5$\\
			\noalign{\smallskip}\noalign{\smallskip}
			$c_4(6,l)$ & $\frac{1}{\sqrt{6}}$ & $\frac{1}{\sqrt{3}}$ & $\frac{1}{\sqrt{2}}$ & $\frac{1}{\sqrt{3}}$ & $\frac{1}{\sqrt{6}}$ & $1$\\
			\noalign{\smallskip}\hline
		\end{tabular}
		\caption{The values of $c_3(6,l)$ and $c_4(6,l)$ for $1\leq l\leq 6$.}
		\label{tab:1}       
	\end{table}
For notational convenience, 
set $a(x)=\frac{e^x}{(1-e^x)^2}$. 
By \eqref{ieq-Ek} we obtain	
\begin{align*}
		|E_6(n)|
		&\le\frac{24N^{2}}{(24n+5)\pi}\exp\left(\frac{\left(24n+5\right)\pi}{12N^2} \right)
		\left(2e^{\frac{5\pi}{18}} + 2e^{\frac{\pi}{9}} \right)
		+2\exp\left(\frac{\pi\left(24n+5\right)}{12N^2} \right)\times \Upsilon
	\end{align*}
	where
	\begin{align*}
		\Upsilon
		&=\frac{2}{\sqrt{6}} \exp\bigg(a(-\pi) +a(-\pi/6) + \frac{5\pi}{144}\bigg)  -  \frac{2}{\sqrt{6}}\exp\left(\frac{5\pi }{144}\right) \nonumber \\[5pt]
		&+\frac{2}{\sqrt{3}}	\exp\bigg(a(-\pi) + a(-2\pi/3) + \frac{\pi}{72}\bigg) -  \frac{2}{\sqrt{3}}\exp\left(\frac{\pi }{72}\right)
		\nonumber \\[5pt]
		&+\frac{1}{\sqrt{2}}
		\exp\bigg(a(-\pi)+ a(-3\pi/2)-\frac{\pi}{48}\bigg)+
		\exp\bigg(a(-\pi)+ a(-6\pi)- \frac{5\pi}{24}\bigg).
	\end{align*}
Now let  $y=\mu_6(n)$, implying that  ${24n+5}=\frac{6^3y^2}{5\pi^5}$. 
By taking $N=\left\lfloor y\right\rfloor$, one can verify that 
	\begin{align*}
		|E_6(n)|
		\leq \frac{5\pi\left\lfloor y\right\rfloor^2}{9y^2}
		\exp\left(\frac{18y^2}{5\pi\left\lfloor y\right\rfloor^2}\right) \left(2e^{\frac{5\pi}{18}} + 2e^{\frac{\pi}{9}} \right)
		+2\exp\left(\frac{18y^2}{5\pi\left\lfloor y\right\rfloor^2}\right)\times\Upsilon.
	\end{align*}
	Then, using the following two inequalities:
	\[\frac{\left\lfloor y\right\rfloor^2}{y^2}\le 1  \quad \text{and} \quad
	\frac{y^2}{\left\lfloor y\right\rfloor^2}<\frac{y^2}{(y-1)^2}< 2\, \quad\text{for }   y\ge 4,\]
	we deduce  that
	\begin{align*}
		|E_6(n)|&\le\frac{10\pi}{9}e^{\frac{36}{5\pi}}\left(e^{\frac{5\pi}{18}} + e^{\frac{\pi}{9}} \right)+2e^{\frac{36}{5\pi}}\times\Upsilon\leq 812
	\end{align*}
holds when $n\ge 3$ satisfies $y=\mu_6(n)\geq 4$.
	
Next we focus on estimating  the value of $|B_6(n)|$. From Proposition \ref{prop-k-n-24} it follows that
	\begin{align*}
		|B_6(n)|\leq	\frac{5\pi^2}{18\sqrt{3}y}
		\sum_{l\in \{1,2,4,5\}}
		\sum_{\substack{{2\le t\le \left\lfloor y\right\rfloor}\\t\equiv_6 l }}I_1\left(\frac{y}{t}\right)
		= \frac{5\pi^2}{18\sqrt{3}y}\sum_{\substack{2\le t\le \left\lfloor y\right\rfloor \\3\nmid t }}I_1\left(\frac{y}{t}\right).
	\end{align*}	
	Thus,
	\begin{align*}
		|B_6(n)|
		\leq \frac{5\pi^2}{18\sqrt{3}y}\frac{2\left\lfloor y\right\rfloor}{3}I_1\left(\frac{y}{2}\right)
		\leq \frac{5\pi^2}{27\sqrt{3}}I_1\left(\frac{y}{2}\right)		\leq \frac{10\pi^{\frac{3}{2}}e^{\frac{y}{2}}}{27\sqrt{3y}},
	\end{align*}
where the last inequality is obtaind by using \eqref{eq-lem-B-K}.
Combining the bounds of $E_6(n)$ and $B_6(n)$, we get
	\begin{align*}
		|R_6(n)| =|E_6(n)| + |B_6(n)|\leq 812 + \frac{10\pi^{\frac{3}{2}}e^{\frac{y}{2}}}{27\sqrt{3y}}
	\end{align*}
whenever $n\ge 3$ satisfies $y=\mu_6(n)\geq 4$. Keep in mind that $y$ depends on $n$.
Thus, it remains to show for $n\geq 46$,
	\begin{equation*}
		\gamma(n):=\frac{10\pi^{\frac{3}{2}}e^{\frac{y}{2}}}{27\sqrt{3y}}=\frac{10\pi^{\frac{3}{2}}e^{\frac{\mu_6(n)}{2}}}{27\sqrt{3\mu_6(n)}}\geq 812.
	\end{equation*}
By studying the derivative of $\gamma(n)$, one can show that it is increasing on the interval $[0,+\infty)$.  
Thus if $n\geq 46$ then $\gamma(n)\geq \gamma(46)>812$, as desired. 
This completes the proof.
\end{proof}

Finally, we are able to prove Theorem \ref{thm-p_k2}.

\begin{proof}[Proof of Theorem \ref{thm-p_k2}]
As before, we only proof the case of $k=6$. 
According to Theorem \ref{thm-p_k1}, if $n\ge 46$, then 
	\begin{align*}
		p_6(n)=M_6(n)+R_6(n),
	\end{align*}
	where $|R_6(n)|\leq \hat{R}_6(n)=\frac{20\pi^{\frac{3}{2}}e^{\frac{y}{2}}}{27\sqrt{3 }y^{\frac{1}{2}}}$.
	
If we let 
	\begin{equation*}
		G_6(n):=\frac{\hat R_6(n)}{M_6(n)}
		= \frac{8\sqrt{2}}{3\sqrt{\pi}}\cdot \frac{\sqrt{y}e^{\frac{y}{2}}}{I_1(y)},
	\end{equation*}
then 
	\begin{equation*}
		M_6(n)(1-G_6(n))\leq p_6(n)\leq M_6(n)(1+G_6(n)).
	\end{equation*}
Note that \eqref{eq-I1s} allows us to deduce that for $n\geq 124$ satisfies $y\geq 26$,
	\begin{align*}
		G_6(n)
		\leq
		\frac{32y^2e^{-\frac{y}{2}}}{3(2y-1)}.	
	\end{align*}
It remains to show 
	\begin{align*}
		\frac{32y^2e^{-\frac{y}{2}}}{3(2y-1)}\leq \frac{1}{y^6}.
	\end{align*}
Keep in mind that $y$ depends on $n$. Thus it is equivalent to show 
	\begin{align*}
		L(n):={32\mu_6(n)^8}-{3(2\mu_6(n)-1)}e^{\frac{\mu_6(n)}{2}}\leq 0,
	\end{align*}
holds for $n\geq 677$. But this can be verified by showing the derivative of $L(n)$
is negative and $L(677)<0$ on the interval $[677,+\infty)$ with the help of mathematical software. This completes the proof.
\end{proof}

\section*{Acknowledgments}
This work was partially supported by the Fundamental Research Funds for the Central Universities.


\begin{thebibliography}{1}
\bibitem{Branden-2015}
P. Br\"and\'en, Unimodality, log-concavity, real-rootedness and beyond, in: M. B\'ona, (ed.), Handbook of Enumerative Combinatorics. Boca Raton, FL: CRC Press, Discrete Math. Appl. 437--483 (2015).

\bibitem{Briggs-1985}
W.E. Briggs, Zeros and factors of polynomials with positive coefficients and protein-ligand binding,
Rocky Mountain J. Math. 15: 75--89 (1985).

\bibitem{Bringmann-Kane-Rolen-Tripp-2021} K. Bringmann, B. Kane, L. Rolen and Z. Tripp, Fractional partitions and conjectures of Chern-Fu-Tang and Heim-Neuhauser, Trans. Amer. Math. Soc. Ser. B 8: 615--634 (2021).

\bibitem{Chern-2019}
S. Chern, Asymptotics for the Fourier coefficients of eta-quotients, J. Number
Theory 199: 168--191 (2019).

\bibitem{Chen-Jia-Wang-2019}
W.Y.C. Chen, D.X.Q. Jia and L.X.W. Wang,
Highter order Tur\'an inequalities for the partition function,
Trans. Amer. Math. Soc. 372(3): 2143--2165 (2019).

\bibitem{Craig-Pun-2021}
W. Craig, A. Pun, A note on the third-order Turán inequalities for $k$-regular partitions, Res. Number
Theory 7(5) (2021).

\bibitem{Corteel-Lovejoy-2004}
S. Corteel, J. Lovejoy, Overpartitions, Trans. Amer. Math. Soc. 356: 1623--1635 (2003).

\bibitem{Csordas-2003}
G. Csordas, Complex zero decreasing sequences and the Riemann hypothesis II, in: Analysis and
Applications—ISAAC 2001 (Berlin), in: Int. Soc. Anal. Appl. Comput., vol. 10, Kluwer Acad. Publ.,
Dordrecht, 121--134 (2003).


\bibitem{DeSalvo-Pak-2015}
S. DeSalvo and I. Pak,
Log-concavity of the partition function,
Ramanujan J. 38(1): 61--73 (2015).

\bibitem{Dimitrov-1998}
D.K. Dimitrov, Higher order Tur\'an inequalities, Proc. Amer. Math. Soc. 126(7): 2033--2037 (1998).

\bibitem{Dong-Ji-2024}
J.J.W. Dong and K.Q. Ji,
Higher order Tur\'an inequalities for the distinct partition function,
J. Number Theory 260: 71--102 (2024).


\bibitem{Engel-2017}
B. Engel, Log-concavity of the overpartition function, Ramanujan J. 43(2): 229--241 (2017).


\bibitem{Fan-Wang}
J.-T. Fan and L.X.W. Wang, Complete monotonicity, $(m,n)$-Laguerre
inequality, Briggs inequality and
Laguerre-P\'olya class, preprint.

\bibitem{G-O-R-Z-2019}
M. Griffin, K. Ono, L. Rolen, and D. Zagier, Jensen polynomials for the Riemann zeta function and other sequences, Proc. Natl. Acad. Sci. USA 116(23): 11103--11110 (2019).

\bibitem{Hou-Zhang-2019}
Q.-H. Hou and Z.-R. Zhang, $r$-log-concavity of partition functions, Ramanujan J. 48
(1): 117--129 (2019).


\bibitem{Jia-Wang-2020}
D.X.Q. Jia and L.X.W. Wang,
Determinantal inequalities for the partition function,
Proc. R. Soc. Edinb. A 150(3): 1451--1466 (2020).

\bibitem{Liu-Zhang-2021} E.Y.S. Liu and H.W.J. Zhang, Inequalities for the overpartition function, Ramanujan J. 54: 485--509 (2021).

\bibitem{Mukherjee-2023}
 G. Mukherjee, Inequalities for the overpartition function arising from determinants, Adv. in Appl.
Math. 152: 102598 (2024).

\bibitem{Mukherjee-Zhang-Zhong-2023}
G. Mukherjee, H.W.J. Zhang and Y. Zhong, Higher order log-concavity of the overpartition function and its consequences, Proc. Edinb. Math. Soc. 66(1): 164--181 (2023).

\bibitem{Niculescu-2000}
C.P. Niculescu, A new look at Newton's inequalities, J. Inequal. Pure Appl. Math. 1(2): 17 (2000).

\bibitem{Nicolas-1978}
J.L. Nicolas, Sur les entiers \textit{N} pour lesquels il y a beaucoup de groupes ab\'eliens d'ordre \textit{N}, Ann. Inst. Fourier 28(4): 1--16 (1978).

\bibitem{Peng-Zhang-Zhong-2023}
Y. Peng, H.W.J. Zhang and Y. Zhong, Inequalities for the $k$-regular overpartitions, arxiv:2308.04678.

\bibitem{Rosset-1989}
S. Rosset, Normalized symmetric functions, Newton's inequalities and a new set of stronger inequalities, Amer. Math. Monthly 96(9): 815--819 (1989).


\bibitem{Wang-Yang-2023}
L.X.W. Wang and N.N.Y. Yang, Positivity of the determinants of the partition function and the overpartition function,
Math. Comp. 92: 1383--1402 (2023).

\bibitem{Wang-Yang}
 L.X.W. Wang and E.Y.Y. Yang, Laguerre inequality and determinantal
inequality for the distinct partition function, submitted.

\bibitem{Zhang-Zhao-2024}
Z.-X. Zhang and J.J.Y. Zhao,
The Briggs inequality of Boros-Moll sequences, arxiv:2402.11620.
\end{thebibliography}
\end{document}